\documentclass[12pt]{article}
\usepackage{amsmath,amsfonts,amssymb,amsthm}
\input amssym.def
\begin{document}
\def \Z{\Bbb Z}
\def \C{\Bbb C}
\def \R{\Bbb R}
\def \Q{\Bbb Q}
\def \N{\Bbb N}
\def \D{\mathcal{D}}
\def \wt{{\rm wt}}
\def \tr{{\rm tr}}
\def \span{{\rm span}}
\def \Res{{\rm Res}}
\def \End{{\rm End}}
\def \Ind {{\rm Ind}}
\def \Irr {{\rm Irr}}
\def \Aut{{\rm Aut}}
\def \Hom{{\rm Hom}}
\def \mod{{\rm mod}}
\def \ann{{\rm Ann}}
\def \<{\langle}
\def \>{\rangle}
\def \g{{\frak{g}}}

\def \be{\begin{equation}\label}
\def \ee{\end{equation}}
\def \bl{\begin{lem}\label}
\def \el{\end{lem}}
\def \bt{\begin{thm}\label}
\def \et{\end{thm}}
\def \bp{\begin{prop}\label}
\def \ep{\end{prop}}
\def \br{\begin{rem}\label}
\def \er{\end{rem}}
\def \bc{\begin{coro}\label}
\def \ec{\end{coro}}
\def \bd{\begin{de}\label}
\def \ed{\end{de}}
\def \bex{\begin{exa}\label}
\def \eex{\end{exa}}

\newtheorem{thm}{Theorem}[section]
\newtheorem{prop}[thm]{Proposition}
\newtheorem{coro}[thm]{Corollary}
\newtheorem{conj}[thm]{Conjecture}
\newtheorem{exa}[thm]{Example}
\newtheorem{lem}[thm]{Lemma}
\newtheorem{rem}[thm]{Remark}
\newtheorem{de}[thm]{Definition}
\newtheorem{hy}[thm]{Hypothesis}
\makeatletter
\@addtoreset{equation}{section}
\def\theequation{\thesection.\arabic{equation}}
\makeatother
\makeatletter

\vspace{5cm}
\begin{center}{\Large \bf Twisted modules and pseudo-endomorphisms}

\vspace{0.5cm}
Haisheng Li\\
Department of Mathematical Sciences\\
Rutgers University, Camden, NJ 08102
\end{center}

\begin{abstract}
We exhibit a connection between two constructions of twisted modules
for a general vertex operator algebra with respect to inner
automorphisms. We also study pseudo-derivations,
pseudo-endomorphisms, and twist deformations of ordinary modules by
pseudo-endomorphisms, which are intrinsically connected to one of
the two constructions.
\end{abstract}

\section{Introduction}
In vertex operator algebra theory, for a vertex operator algebra
$V$, in addition to the notion of $V$-module there is a notion of
$\sigma$-twisted $V$-module where $\sigma$ is a finite order
automorphism of $V$. This is one of those new features of vertex
operator algebras, in contrast with classical Lie or associative
algebras. The notion of twisted module was originated from the
construction of the celebrated moonshine module vertex operator
algebra $V^{\natural}$ (see \cite{flm1}, \cite{flm}, \cite{le}) and
it had played a key role therein. Since then, the theory of twisted
representations has been extensively studied in literature (cf.
[FFR], [D], \cite{li-twisted}, \cite{dlm-twisted}).

In \cite{li-twisted}, a conceptual construction of vertex algebras
and their twisted modules was obtained and  a canonical construction
of twisted modules with respect to inner automorphisms was found.
Let $V$ be a general vertex operator algebra and let $h\in V$,
satisfying that $L(n)h=\delta_{n,0}h$ and $h(n)h=\delta_{n,1}\alpha$
for $n\in \N$, where $\alpha$ is a rational number. Assume that
$h(0)$ acts semisimply on $V$ with eigenvalues in $\frac{1}{r}\Z$
for some positive integer $r$. Note that $e^{2\pi i h(0)}$, denoted
by $\sigma_{h}$, is an automorphism of $V$. Set
$$\Delta(h,z)=z^{h(0)}\exp\left(\sum_{n\ge 1}\frac{h(n)}{-n}z^{-n}\right).$$
It was proved therein that for any $V$-module $(W,Y_{W}(\cdot,z))$,
$(W,Y_{W}(\Delta(h,z)\cdot,z))$ carries the structure of a
$\sigma_{h}$-twisted $V$-module. This correspondence gives rise to a
canonical isomorphism between the category of $V$-modules and the
category of $\sigma_{h}$-twisted $V$-modules.

There is another construction of twisted modules by employing the
construction of contragredient modules. Let $V$ be a vertex operator
algebra. A fundamental result due to \cite{fhl} is that for any
$V$-module $W$, the restricted dual $W'$ of $W$ has a natural
$V$-module structure, called the contragredient module. Note that
the construction of contragredient module depends on the conformal
vector; it depends {\em not only} on the vertex algebra structure.
In \cite{dlm-sl2}, to study modules for vertex operator algebras
associated to affine Lie algebra $\hat{sl}_{2}$ at admissible
levels, we employed a different conformal vector whose $L(0)$-weight
grading corresponds to the principal grading, whereas the usual
Segal-Sugawara conformal vector corresponds to the homogeneous
grading. More generally, let $V$ be a vertex operator algebra with
conformal vector $\omega$ and let $h\in V$ be such that
$L(n)h=\delta_{n,0}h$ and $h(n)h=\delta_{n,1}\alpha$ for $n\in \N$,
where $\alpha$ is a rational number. Then $\omega+L(-1)h$ is a
conformal vector. (This is a reformulation of the Feigin-Fuchs
construction of modules for the Virasoro algebra.) In general,
vertex algebra $V$ equipped with the new conformal vector
$\omega+L(-1)h$ is a $\Q$-graded vertex operator algebra, whose
notion differs from the ordinary one only on the grading group. It
was observed therein that for a general $\Q$-graded vertex operator
algebra $V$, essentially the same argument of \cite{fhl} shown that
the contragredient dual of a $V$-module is a twisted module with
respect to the automorphism $\sigma=e^{4\pi iL(0)}$. Then, given an
ordinary vertex operator algebra $V$ with an element $h\in V$
satisfying the conditions as above, we have a different construction
of twisted modules by considering contragredient dual of $V$-modules
with respect to the new conformal vector associated to $h$.

This present paper is to study the relation between the two
constructions of twisted modules and to study a certain general
theory closely related to the first construction. More specifically,
in the first part, we give a connection between the two
constructions of twisted modules and in the second we study
pseudo-endomorphisms and the twist deformations of ordinary modules
by pseudo-endomorphisms.

Note that the first construction is simply to twist ordinary
representations.  Closely related to this spirit is the physics
super-selection theory (cf. \cite{hk}, \cite{frs}).  For an
illustration, let $A$ be an associative ($C^{*}$-)algebra. A trivial
fact is that for any endomorphism $\sigma$ of $A$ and for any
representation $\rho$ of $A$, $\rho\circ \sigma$ is also a
representation of $A$. In general, $A$ has infinitely many
nonequivalent irreducible modules, but in practice only certain
irreducible modules are selected by certain rules. The super
selection rule is to select those representations obtained from the
distinguished ``vacuum module'' through (unitary) endomorphisms.
Also in this very theory, fusion rule can be defined and studied in
terms of endomorphisms.  For vertex operator algebras, the adjoint
modules are understood to be the vacuum modules. All of these
motivated the work \cite{li-selection} in which the physics
super-selection theory was studied in the contexts of vertex
operator algebras. The key idea is to replace endomorphisms with
``pseudo-endomorphisms'' like $\Delta(h,z)$, which are not
endomorphisms in the usual sense. (Notice that if $V$ is a simple
vertex operator algebra, the $\sigma$-twist of $V$ is always
isomorphic to $V$ for any endomorphism $\sigma$ of $V$.) The
``pseudo-endomorphisms'' considered in \cite{li-selection} are
elements $\Delta(z)\in \Hom (V,V\otimes \C[z,z^{-1}])$, satisfying
\begin{eqnarray*}
&&\Delta(z){\bf 1}={\bf 1},\ \ \ \
[L(-1),\Delta(z)]=-\frac{d}{dz}\Delta(z),\\
&&\Delta(z)Y(v,z_{0})=Y(\Delta(z+z_{0})v,z_{0})\Delta(z)
\end{eqnarray*}
for $v\in V$. It was proved therein that for any $V$-module
$(W,Y_{W}(\cdot,z))$, the pair $(W,Y_{W}(\Delta(z)\cdot,z))$ also
carries the structure of a $V$-module. Furthermore, fusion rule in
the sense of \cite{fhl} was studied in terms of
``pseudo-endomorphisms.'' In particular, it was proved that twisting
the adjoint module $V$ by  $\Delta(h,z)$ gives simple currents in a
certain sense.

The notion of pseudo-endomorphism was formally introduced and
further studied later in \cite{li-pseudo}. Motivated by
Etingof-Kazhdan's notion of pseudo-derivation (see \cite{ek}), we
defined a pseudo-endomorphism  of $V$ more generally as an element
$\Delta(z)\in \Hom (V,V\otimes \C((z)))$ with the same set of
axioms, whereas a pseudo-derivation of $V$ is an element $\psi(z)\in
\Hom (V,V\otimes \C((z)))$ such that
$$[L(-1),\psi(z)]=-\frac{d}{dz}\psi(z),\ \ \ \
[\psi(z),Y(v,x)]=Y(\psi(z+x)v,x)$$ for $v\in V$. Just as with the
classical case, the exponential of a pseudo-derivation, provided it
exists, is a pseudo-automorphism. A result proved in \cite{ek} is
that for any $f(z)\in \C((z)),\ v\in V$,
$$X_{v,f}(z):=\sum_{n\ge 0}\frac{f^{(n)}(z)}{n!}v_{n}$$
is a pseudo-derivation of $V$.  As it was observed in
\cite{li-pseudo} (Remark 3.9), if we could take $f(z)=\log z$, we
would have
$$X_{h,f}=h_{0}\log z+\sum_{n\ge 1}(-1)^{n-1}\frac{1}{n}h_{n}z^{-n}$$
and $\Delta(h,z)=e^{X_{h,f}(z)}$. As the exponentials of
pseudo-derivations ``are'' pseudo-automorphisms, the importance of
exponentiating objects like $X_{h,f}$ is manifest. This suggests to
study ``logarithmic'' pseudo-derivations and pseudo-endomorphisms,
which are elements of $\Hom (V,V\otimes
 \C((z))[\log z])$.

The second part of this present paper is a continuation of
\cite{li-pseudo}, to study pseudo-derivations, their formal
exponentials, and pseudo-endomorphisms. More specifically, we study
more general pseudo-derivations and pseudo-endomorphisms with
$\C((z))$ replaced by certain commutative algebras $R$ over
$\C((z))$ such as $\C((z))[\log z]$ and $\C((z^{1/r}))$ with $r$ a
positive integer. We also study pseudo-endomorphism twists of
ordinary representations for vertex algebras. For the vertex
operator algebra associated to a Heisenberg algebra, we show that
certain pseudo-endomorphism twists of the adjoint module gives
non-highest weight modules.

Recently,  Huang (see \cite{huang}) obtained a very interesting
generalization of the first construction, where $h(0)$ is allowed to
be locally finite, instead of semisimple. In his generalization, the
factor $z^{h(0)}$ in the construction is essentially replaced with
$e^{h(0)\log z}$ in a certain natural way. Then Huang obtained what
he called generalized twisted modules with respect to an
automorphism which is no longer of finite order.

This paper is organized as follows:  In Section 2, we give a
connection between the two constructions of twisted modules. In
Section 3, we study more general pseudo-derivations and
pseudo-endomorphisms.

\section{Two constructions of $\sigma_{h}$-twisted modules}

In this section,  we first review the two constructions of twisted
modules for an inner automorphism and we then give a connection
between the two constructions.

We start with some basic notions. An {\em automorphism} of a vertex
algebra $V$ is a bijective linear endomorphism $\sigma$ of $V$ such
that $\sigma({\bf 1})={\bf 1}$ and
$\sigma(Y(u,z)v)=Y(\sigma(u),z)\sigma(v)$ for $u,v\in V$. An {\em
automorphism} of a vertex operator algebra $V$ is an automorphism
$\sigma$ of $V$ viewed as a vertex algebra such that
$\sigma(\omega)=\omega$, where $\omega$ is the conformal vector of
$V$. Consequently, every automorphism of a vertex operator algebra
preserves all the homogeneous ($L(0)$-weight) subspaces.

A simple fact is that for any vertex algebra $V$ and for any $a\in
V$, $e^{a_{0}}$ is an automorphism of $V$, provided that $a_{0}$ is
locally finite on $V$. Assume that $V$ is a vertex operator algebra.
Let $a\in V_{(1)}$, i.e., $a\in V$ with $L(0)a=a$. It follows that
$a_{0}$ preserves the $L(0)$-weight spaces of $V$, which implies
that $a_{0}$ is locally finite on $V$. Thus $e^{a_{0}}$ is an
automorphism of $V$ viewed as a vertex algebra. Furthermore, if $a$
is primary in the sense that $L(n)a=0$ for $n\ge 1$, we have
$e^{a_{0}}(\omega)=\omega$, so that $e^{a_{0}}$ is an automorphism
of vertex operator algebra $V$.

Let $V$ be a vertex algebra and let $\sigma$ be an automorphism of
$V$ of order $k$. A {\em $\sigma$-twisted $V$-module} (see
\cite{flm}, \cite{ffr}, \cite{dong}) is a vector space $W$, equipped
with a linear map
\begin{eqnarray*}
Y_{W}:&& V\rightarrow \Hom (W,W((x^{1/k})))\subset (\End
W)[[x^{1/k}, x^{-1/k}]]\\
&&v\mapsto Y_{W}(v,x),
\end{eqnarray*}
satisfying the conditions
that $Y_{W}({\bf 1},x)=1$,
$$Y_{W}(\sigma v,x)=\lim_{x^{1/k}\rightarrow \omega_{k}^{-1}x^{1/k}}Y_{W}(v,x)$$
for $v\in V$, where $\omega_{k}=\exp (2\pi i/k)$, and that for $u\in
V^j, v\in V$ with $0\le j\le k-1$,
\begin{eqnarray}\label{jacobi}
&&x_{0}^{-1}\delta\left(\frac{x_{1}-x_{2}}{x_{0}}\right)
Y_W(u,x_{1})Y_W(v,x_{2})-x_{0}^{-1}\delta\left(\frac{x_{2}-x_{1}}{-x_{0}}\right)
Y_W(v,x_{2})Y_W(u,x_{1})\nonumber\\
&&\ \ \ \ =x_2^{-1}\left(\frac{x_1-x_0}{x_2}\right)^{-\frac{j}{k}}
\delta\left(\frac{x_1-x_0}{x_2}\right) Y_W(Y(u,x_0)v,x_2)
\end{eqnarray}
({\em $\sigma$-twisted Jacobi identity}), where $$V^j=\{a\in V\;|\;
\sigma(a)=\omega_{k}^{j}a\}.$$ If $V$ is a vertex operator algebra,
a $\sigma$-twisted $V$-module, which by definition is $\C$-graded by
the $L(0)$-weight, satisfies the two grading restrictions.

Let $V$ be a vertex operator algebra and let $h\in V$ be such that
\begin{eqnarray}\label{e2.15}
L(n)h=\delta_{n,0}h,\ \ h_{n}h=\delta_{n,1}\gamma {\bf 1}\ \
\mbox{for  }n\in \N,
\end{eqnarray}
where $\gamma$ is a fixed rational number, and such that $h_{0}$
acts semisimply on $V$ with rational eigenvalues. In view of
Borcherds' commutator formula we have
\begin{eqnarray}
[L(m),h_{n}]=-nh_{m+n},\ \ \ \ [h_{m},h_{n}]=m\gamma \delta_{m+n,0}
\end{eqnarray}
for $m,n\in \Z$. We shall also freely use $h(n)$ for $h_{n}$.

Set
\begin{eqnarray}
\Delta(h,z)=z^{h(0)}\exp\left(\sum_{k=1}^{\infty}\frac{h(k)}{-k}
(-z)^{-k}\right)\in ({\rm End}V)\{z\}.
\end{eqnarray}
Note that $e^{2\pi ih(0)}$ is an automorphism of $V$. Set
$\sigma_{h}=e^{2\pi ih(0)}$ and assume that $\sigma_{h}$ is of
finite order. The following proposition was proved in
\cite{li-twisted}:

\bp{p3.1} Let $V$ be a vertex operator algebra and let $h\in V$ be
such that (\ref{e2.15}) holds and such that $h(0)$ acts semisimply
on $V$ with rational eigenvalues. Let $(W,Y_{W}(\cdot,z))$ be a
$V$-module. Then $(W, Y_{W}(\Delta(h,z)\cdot,z))$ carries the
structure of a $\sigma_{h}$-twisted $V$-module. \ep

The following is essentially the Feigin-Fuchs bosonization of
Virasoro algebra (cf. \cite{dlinm}):

\bp{p3.4} Let $V$ be a vertex operator algebra with conformal vector
$\omega$ of central charge $c$ and $h\in V$, satisfying
(\ref{e2.15}). Then $\omega+{1\over 2}L(-1)h$ is a conformal vector
of central charge $c-3\gamma$. \ep

In Proposition \ref{p3.4},  $V$ with the new conformal vector
$\tilde{\omega}$ in general is a ${\Q}$-graded vertex operator
algebra in the following sense (see \cite{dlm-sl2}):

\bd{dqvoa} {\em If $V$ satisfies all the axioms of a vertex operator
algebra except that $V$ is $\Q$-graded instead of $\Z$-graded, we
call $V$ {\em a $\Q$-graded vertex operator algebra}.} \ed

Let $V$ be a general ${\Q}$-graded vertex operator algebra. For
$u\in V_{(\alpha)}$ with $\alpha\in \Q$ and for $n\in \Z$, we have
(see \cite{fhl})
\begin{eqnarray}
[L(0),u_{n}]=(\alpha-n-1)u_{n}.
\end{eqnarray}
For $u\in V_{(\alpha)},\; v\in V_{(\beta)}$ with $\alpha,\beta \in
\Q$, we have
\begin{eqnarray*}
e^{2\pi iL(0)}(u_{n}v)=e^{2(\alpha+\beta-n-1)\pi i}(u_{n}v)
=e^{2(\alpha+\beta)\pi i}(u_{n}v) =(e^{2\pi iL(0)}u)_{n}(e^{2\pi
iL(0)}v).
\end{eqnarray*}
Since $L(0){\bf 1}=\omega_{1}{\bf 1}=0$ and $L(0)\omega=2\omega$, we
have
$$e^{2\pi iL(0)}({\bf 1})={\bf 1},\ \ e^{2\pi iL(0)}(\omega)=\omega.$$
Thus $\tau:=e^{2\pi i L(0)}$ is an automorphism of $V$. It is clear
that $\tau$ is of finite order if and only if there is a positive
integer $T$ such that all the $L(0)$-weights, i.e., the
$L(0)$-eigenvalues, are contained in ${1\over T}{\Z}$, which is true
if $V$ as a vertex algebra is finitely generated.

Now, we assume that $V=\oplus_{\alpha\in
\frac{1}{2r}\Z}V_{(\alpha)}$ for some positive integer $r$. Then
$$\tau^{2r}=e^{4r\pi iL(0)}=1\   \mbox{ on }V.$$
That is, $2r$ is a period of $\tau$ and $r$ is a period of
$\tau^{2}$.

Let $W=\oplus_{\alpha\in \C}W_{(\alpha)}$ be a $V$-module. Following
\cite{fhl}, set
$$W'=\oplus_{\alpha\in \C}W_{(\alpha)}^{*}$$
(the {\em restricted dual} of $W$) and define
\begin{eqnarray}\label{edefinition}
\<Y_{W'}(v,z)w',w\> = \<w',Y_{W}(e^{zL(1)}e^{\pi i
L(0)}z^{-2L(0)}v,z^{-1})w\>
\end{eqnarray}
for $v\in V,\; w'\in W',\; w\in W$.  {}From the argument in
\cite{fhl}, we also have
\begin{eqnarray}\label{eback-relation}
\<w',Y(v,z)w\> = \<Y'(e^{zL(1)}e^{-\pi i
L(0)}z^{-2L(0)}v,z^{-1})w',w\>.
\end{eqnarray}

The following was observed in \cite{dlm-sl2}:

\bp{p3.2} Let $V$ be a $\Q$-graded vertex operator algebra such that
$V=\oplus_{\alpha\in \frac{1}{2r}\Z}V_{(\alpha)}$ for some positive
integer $r$ and let $W$ be a $V$-module. Then $(W', Y_{W'})$ carries
the structure of a $\sigma$-twisted $V$-module with $\sigma=e^{4\pi
iL(0)}\;(=\tau^{2})$. \ep

\begin{proof} It is just a slight modification of the proof of the
contragredient module theorem in \cite{fhl}. First, as the
$L(0)$-eigenvalues are contained in $\frac{1}{2r}\Z$, we see {}from
(\ref{edefinition}) that for $v\in V$, $$Y_{W'}(v,z)\in (\End
W)[[z^{1/r},z^{-1/r}]]$$ and we have
$$Y_{W'}(v,z)\in \Hom(W,W((z^{1/r})))$$
{}from the grading restrictions on $W$ as in \cite{fhl}. Second, for
$v\in V$ we have
$$\lim_{z^{1/r}\rightarrow \omega_{r}^{-1}z^{1/r}}z^{-2L(0)}v
=\lim_{z^{1/r}\rightarrow
\omega_{r}^{-1}z^{1/r}}\left(z^{1/r}\right)^{-2rL(0)}v=e^{4\pi
iL(0)}z^{-2L(0)}v,$$ so that
$$Y_{W'}(\sigma v,z)=Y_{W'}(e^{4\pi iL(0)}v,z)
=\lim_{z^{1/r}\rightarrow \omega_{r}^{-1}z^{1/r}} Y_{W'}(v,z),$$
where $\omega_{r}=\exp (2\pi i/r)\in \C$. Then it remains to prove
the twisted Jacobi identity. Let $a\in V_{(\alpha)},\; b\in V$ with
$\alpha\in \frac{1}{2r}\Z$. We need to prove
\begin{eqnarray}\label{jacobi-proof}
&&z_{0}^{-1}\delta\left(\frac{z_{1}-z_{2}}{z_{0}}\right)Y_{W'}(a,z_{1})Y_{W'}(b,z_{2})
-z_{0}^{-1}\delta\left(\frac{z_{2}-z_{1}}{-z_{0}}\right)
Y_{W'}(b,z_{2})Y_{W'}(a,z_{1})\nonumber\\
&&\ \ =z_{2}^{-1}\left(\frac{z_{1}-z_{0}}{z_{2}}\right)^{-2\alpha}
\delta\left(\frac{z_{1}-z_{0}}{z_{2}}\right)
Y_{W'}(Y(a,z_{0})b,z_{2}).
\end{eqnarray}
 Recall the following formulas
\begin{eqnarray}
z^{L(0)}L(-1)&=&z^{-1}L(-1)z^{L(0)},\label{3.3}\\
e^{zL(1)}L(-1)&=&L(-1)e^{zL(1)}+2ze^{zL(1)}L(0)-z^{2}
e^{zL(1)}L(1),\\
z^{L(0)}Y(a,z_{0})z^{-L(0)}&=&Y(z^{L(0)}a,z_{0}z),\\
e^{zL(1)}Y(a,z_{0})e^{-zL(1)}&=&Y\left(e^{z(1-zz_{0})}
(1-z_{0}z)^{-2L(0)}a, \frac{z_{0}}{1-zz_{0}}\right).\label{3.6}
\end{eqnarray}

For $u,v\in V$, from \cite{fhl} we have
\begin{eqnarray}
& &z_{0}^{-1}\delta\left(\frac{z_{1}-z_{2}}{z_{0}}\right)
Y_{W}\left(v,z_{2}^{-1}\right) Y_{W}\left(u,z_{1}^{-1}\right)
\nonumber\\
& &\ \ \ \ -z_{0}^{-1}\delta\left(\frac{z_{2}-z_{1}}{-z_{0}}\right)
Y_{W}\left(u,z_{1}^{-1}\right) Y_{W}\left(v,z_{2}^{-1}\right)
\nonumber\\
&=&-(z_{1}z_{2})^{-1}\left(-\frac{z_{0}}{z_{1}z_{2}}\right)^{-1}\delta\left(
\frac{-z_{2}^{-1}+z_{1}^{-1}}{-{z_{0}\over z_{1}z_{2}}}\right)
Y_{W}\left(v,z_{2}^{-1}\right)Y_{W}\left(u,z_{1}^{-1}\right)\nonumber\\
& &+(z_{1}z_{2})^{-1}\left(-\frac{z_{0}}{z_{1}z_{2}}\right)^{-1}\delta\left(
\frac{z_{1}^{-1}-z_{2}^{-1}}{-{z_{0}\over z_{1}z_{2}}}\right)
Y_{W}\left(u,z_{1}^{-1}\right)Y_{W}\left(v,z_{2}^{-1}\right)\nonumber\\
&=&(z_{1}z_{2})^{-1}z_{2}\delta\left(\frac{z_{2}^{-1}-{z_{0}\over
z_{1}z_{2}}} {z_{1}^{-1}}\right)Y_{W}\left(Y\left(u,
-\frac{z_{0}}{z_{1}z_{2}}\right)v,z_{2}^{-1}\right)\nonumber\\
&=&z_{1}^{-1}\delta\left(\frac{z_{2}+z_{0}}{z_{1}}\right)
Y_{W}\left(Y\left(u,-\frac{z_{0}}{z_{1}z_{2}}\right)v,z_{2}^{-1}
\right).
\end{eqnarray}
Now, take
\begin{eqnarray*}
u=e^{z_{1}L(1)}e^{\pi iL(0)}z_{1}^{-2L(0)}a,\ \
v=e^{z_{2}L(1)}e^{\pi iL(0)}z_{2}^{-2L(0)}b.
\end{eqnarray*}
For proving the desired twisted Jacobi identity, it suffices to
prove
\begin{eqnarray*}
& &z_{2}^{-1}\delta\left(\frac{z_{1}-z_{0}}{z_{2}}\right)
\left(\frac{z_{1}-z_{0}}{z_{2}}\right)^{-2\alpha}
e^{z_{2}L(1)}e^{\pi iL(0)}z_{2}^{-2L(0)}Y(a,z_{0})\nonumber\\
&=&z_{1}^{-1}\delta\left(\frac{z_{2}+z_{0}}{z_{1}}\right)
Y\left(e^{z_{1}L(1)}e^{\pi iL(0)}z_{1}^{-2L(0)}a,
-\frac{z_{0}}{z_{1}z_{2}}\right)e^{z_{2}L(1)}e^{\pi iL(0)}
z_{2}^{-2L(0)},
\end{eqnarray*}
or, equivalently to prove
\begin{eqnarray}
& &z_{2}^{-1}\delta\left(\frac{z_{1}-z_{0}}{z_{2}}\right)
\left(\frac{z_{1}-z_{0}}{z_{2}}\right)^{-2\alpha}\cdot\nonumber\\
& &\ \ \ \ \cdot e^{z_{2}L(1)}e^{\pi iL(0)}z_{2}^{-2L(0)}
Y(a,z_{0})z_{2}^{2L(0)}
e^{-\pi iL(0)}e^{-z_{2}L(1)}\nonumber\\
&=&z_{1}^{-1}\delta\left(\frac{z_{2}+z_{0}}{z_{1}}\right)
Y\left(e^{z_{1}L(1)}e^{\pi iL(0)}z_{1}^{-2L(0)}a,
-\frac{z_{0}}{z_{1}z_{2}}\right).
\end{eqnarray}
Using formulas (\ref{3.3})-(\ref{3.6}), we get
\begin{eqnarray}
& &e^{z_{2}L(1)}e^{\pi iL(0)}z_{2}^{-2L(0)}Y(a,z_{0}) z_{2}^{2L(0)}
e^{-\pi iL(0)}e^{-z_{2}L(1)}\nonumber\\
&=&e^{z_{2}L(1)}e^{\pi iL(0)}Y\left(z_{2}^{-2L(0)}a,
-z_{2}^{-2}z_{0}\right)
e^{-\pi iL(0)}e^{-z_{2}L(1)}\nonumber\\
&=&e^{z_{2}L(1)}Y\left(e^{\pi iL(0)}z_{2}^{-2L(0)}a,
-z_{2}^{-2}z_{0}
\right)e^{-z_{2}L(1)}\nonumber\\
&=&Y\left(e^{z_{2}(1+z_{0}z_{2}^{-1})L(1)}(1+z_{0}z_{2}^{-1})^{-2L(0)}
e^{\pi iL(0)}
z_{2}^{-2L(0)}a,-\frac{z_{2}^{-2}z_{0}}{1+z_{0}z_{2}^{-1}}\right)
\nonumber\\
&=&Y\left(e^{(z_{2}+z_{0})L(1)}e^{\pi iL(0)}
(z_{2}+z_{0})^{-2L(0)}a,-\frac{z_{0}}{z_{2}(z_{2}+z_{0})}\right).
\end{eqnarray}
Thus
\begin{eqnarray}
& &z_{2}^{-1}\delta\left(\frac{z_{1}-z_{0}}{z_{2}}\right)
\left(\frac{z_{1}-z_{0}}{z_{2}}\right)^{-2\alpha}\cdot\nonumber\\
& &\cdot e^{z_{2}L(1)}e^{\pi iL(0)}z_{2}^{-2L(0)}
Y(a,z_{0})z_{2}^{2L(0)}
e^{-\pi iL(0)}e^{-z_{2}L(1)}\nonumber\\
&=&z_{1}^{-1}\delta\left(\frac{z_{2}+z_{0}}{z_{1}}\right)
\left(\frac{z_{2}+z_{0}}{z_{1}}\right)^{2\alpha}\cdot\nonumber\\
& &\cdot Y\left(e^{(z_{2}+z_{0})L(1)}e^{\pi iL(0)}
(z_{2}+z_{0})^{-2L(0)}a,-\frac{z_{0}}{z_{2}(z_{2}+z_{0})}\right)\nonumber\\
&=&z_{1}^{-1}\delta\left(\frac{z_{2}+z_{0}}{z_{1}}\right)
Y\left(e^{z_{1}L(1)}e^{\pi iL(0)}
z_{1}^{-2L(0)}a,-\frac{z_{0}}{z_{2}z_{1}}\right).
\end{eqnarray}
Now the proof is complete.
\end{proof}

\br{rdouble-dual} {\em Let $V$ be a vertex operator algebra and let
$W$ be a $V$-module. It was proved in \cite{fhl} that the double
contragredient module $(W')'$ is isomorphic to $W$. Now, assume that
$V$ is a $\Q$-graded vertex operator algebra. Following the argument
in \cite{fhl}, one obtains
\begin{eqnarray}
(Y')'(v,z)=Y(e^{2\pi iL(0)}v,z)\;\;\left(=Y(\tau (v),z)\right)\ \ \
\mbox{ for }v\in V.
\end{eqnarray}
That is, the double contragredient module $(W')'$ is isomorphic to
the $\tau$-twist of $W$, which is defined by $W^{\tau}=W$ as a
vector space and $Y_{W^{\tau}}(v,z)=Y_{W}(\tau(v),z)$ for $v\in V$.}
\er

The following is the main result of this section, which gives a
connection between the two constructions of $\sigma_{h}$-twisted
modules:

\bt{t3.6} Let $(V,Y, {\bf 1},\omega)$ be a vertex operator algebra
and let $h\in V$, satisfying (\ref{e2.15}). Set
$\tilde{\omega}=\omega+{1\over 2}L(-1)h$. Let $W$ be a $V$-module
and let $(W', Y')$ denote the contragredient module of $W$. Identity
$W$ with $(W')'$ canonically as a vector space. Then
$(W,(Y')_{\tilde{\omega}}')$ is isomorphic to $(W,
Y(\Delta(h,z)e^{-{1\over 2}\pi ih(0)}\cdot,z))$ as a
$\sigma_{h}$-twisted $V$-module, where $(Y')'_{\tilde{\omega}}$
denotes the contragredient dual of $Y'$ with respect to the
conformal vector $\tilde{\omega}$. \et

\begin{proof} First from \cite{fhl} we have
\begin{eqnarray}
z^{L(0)}e^{z_{0}L(1)}=e^{z_{0}z^{-1}L(1)}z^{L(0)}.
\end{eqnarray}
With $\tilde{\omega}=\omega +{1\over 2}L(-1)h$, we have
\begin{eqnarray}
\tilde{L}(m)=L(m)-{1\over 2}(m+1)h(m)\;\;\;\mbox{ for }m\in {\Z},
\end{eqnarray}
where $Y(\tilde{\omega},z)=\sum_{n\in \Z}\tilde{L}(n)z^{-n-2}$. In
particular, we have
\begin{eqnarray}\label{e3.8}
\tilde{L}(-1)=L(-1),\ \ \tilde{L}(0)=L(0)-{1\over 2}h(0),\ \
\tilde{L}(1)=L(1)-h(1).
\end{eqnarray}
Let $a\in V,\; u\in W,\; v'\in W'$. Using (\ref{eback-relation}) we
have
\begin{eqnarray}
& &\langle v',Y(a,z)u\rangle\nonumber\\
&=&\langle Y'(e^{zL(1)}(-z^{-2})^{L(0)}a,z^{-1})v',u\rangle\nonumber\\
&=&\langle v',(Y')_{\tilde{\omega}}'(e^{z^{-1}\tilde{L}(1)}e^{-\pi
i\tilde{L}(0)}z^{2\tilde{L}(0)}
e^{zL(1)}(-z^{-2})^{L(0)}a,z)u\rangle.
\end{eqnarray}
Recall that $[h(0),h(1)]=0$ and $[h(0), L(j)]=0$ for $j=0,1$ and
also recall {}from \cite{fhl} (formula (5.3.1)) that
$$(-z^{2})^{L(0)} e^{zL(1)}(-z^{-2})^{L(0)}=e^{-z^{-1}L(1)}.$$
Then
\begin{eqnarray}
& &e^{z^{-1}\tilde{L}(1)}e^{-\pi i\tilde{L}(0)}z^{2\tilde{L}(0)}
e^{zL(1)}(-z^{-2})^{L(0)}\nonumber\\
&=&e^{z^{-1}(L(1)-h(1))}e^{-\pi i(L(0)-{1\over 2}h(0))}
z^{2(L(0)-{1\over 2}h(0))}
e^{zL(1)}(-z^{-2})^{L(0)}\nonumber\\
&=&z^{-h(0)}e^{z^{-1}(L(1)-h(1))}(-z^{2})^{L(0)}
e^{zL(1)}(-z^{-2})^{L(0)}e^{{1\over 2}\pi ih(0)}\nonumber\\
&=&z^{-h(0)}e^{z^{-1}(L(1)-h(1))}e^{-z^{-1}L(1)}e^{{1\over 2}\pi
ih(0)}.
\end{eqnarray}
By Lemma 3.2 from \cite{li-ext} we have
\begin{eqnarray}
e^{z(L(1)-h(1))}e^{-zL(1)}=\exp\left(\sum_{k=1}^{\infty}
\frac{h(k)}{k}(-z)^{k}\right).
\end{eqnarray}
Using this we obtain
\begin{eqnarray*}
& &\< v',Y(a,z)u\>\nonumber\\
&=&\<
v',(Y')_{\tilde{\omega}}'\left(z^{-h(0)}e^{z^{-1}(L(1)-h(1))}e^{-z^{-1}L(1)}
e^{{1\over 2}\pi ih(0)}a,z\right)u\> \nonumber\\
&=&\<v',(Y')_{\tilde{\omega}}'\left(z^{-h(0)}\exp\left(\sum_{k=1}^{\infty}
\frac{h(k)}{k}(-z^{-1})^{k}\right)e^{{1\over 2}\pi
ih(0)}a,z\right)u\>
\nonumber\\
&=&\<v',(Y')_{\tilde{\omega}}'\left(\Delta(-h,z)e^{{1\over 2}\pi
ih(0)}a,z\right)u\>.
\end{eqnarray*}
Thus we get
\begin{eqnarray}
\< v',Y\left(\Delta(h,z)e^{-{1\over 2}\pi ih(0)}a,z\right)u\> =\<
v',(Y')_{\tilde{\omega}}'(a,z)u\>,
\end{eqnarray}
as desired.
\end{proof}

\section{Pseudo-derivations and pseudo-endomorphisms}
In this section, we study what we call $R$-valued pseudo-derivations
and pseudo-endomorphisms with $R$ a commutative associative algebra
over $\C((z))$. These notions naturally extend those of
($\C((z))$-valued) pseudo-derivation in \cite{ek} and
pseudo-endomorphism in \cite{li-pseudo}. We also study
pseudo-endomorphism twists of ordinary representations for a vertex
algebra.

Let $R$ be a commutative associative algebra over $\C((z))$,
equipped with a $\C$-linear derivation $\partial$ on $R$, satisfying
$$\partial f(z)=\frac{d}{dz}f(z)\ \ \ \mbox{ for }f(z)\in \C((z)).$$
Typical examples of such an $R$ are $\C((z))$, $\C((z^{1/r}))$ with
$r$ a positive integer, and $\C((z^{1/r}))[\log z]$ with $\log z$ as
a new variable where $\partial \log z=z^{-1}$.

We now fix such an $R$ together with $\partial$ throughout this
section.

\bd{d-rpseudo-der} {\em Let $V$ be a vertex algebra. An {\em
$R$-valued pseudo-derivation} of $V$ is an element
$$\Psi\in \Hom (V,V\otimes R),$$ satisfying
\begin{eqnarray}\label{e-pseudo-der}
 \Psi Y(u,x)v-(Y(u,x)\otimes 1)\Psi(v)=Y(e^{x(1\otimes
 \partial)}\Psi(u),x)(v\otimes 1)
\end{eqnarray}
for $u,v\in V$, where for the expression on right it is understood
that
\begin{eqnarray*}
Y(u\otimes a,x)(v\otimes b)=Y(u,x)v\otimes ab
\end{eqnarray*}
for $u,v\in V,\; a,b\in R$.} \ed

\bd{d-rpseudo-end} {\em Let $V$ be a vertex algebra. An {\em
$R$-valued pseudo-endomorphism} of $V$ is an element $\Delta\in \Hom
(V,V\otimes R),$ satisfying
\begin{eqnarray}
\Delta ({\bf 1})={\bf 1}\otimes 1,\ \ \ \ \Delta
Y(v,x)=Y(e^{x(1\otimes \partial)}\Delta (v),x)\Delta
\end{eqnarray}
for $v\in V$.} \ed

As consequences of the definitions we have (cf. \cite{li-pseudo}):

\bl{ldproperty-der} Let $V$ be a vertex algebra. For any $R$-valued
pseudo-derivation $\Psi$ of $V$, we have
\begin{eqnarray}
\Psi ({\bf 1})=0,\ \ \ \ (\D\otimes 1)\Psi-\Psi \D=-(1\otimes
\partial)\Psi,
\end{eqnarray}
where $\D$ is the linear operator on $V$ defined by $\D v=v_{-2}{\bf
1}$ for $v\in V$. \el

\begin{proof} It is basically the same proof as in \cite{li-pseudo} (Proposition 2.3).
We have $\Psi ({\bf 1})=0$ because
$$Y(e^{x(1\otimes \partial)}\Psi({\bf 1}),x)=\Psi Y({\bf
1},x)-(Y({\bf 1},x)\otimes 1)\Psi=\Psi-\Psi=0$$ and because $Y$ is
an injective map. By definition we have
$$\Psi Y(v,x){\bf 1}-(Y(v,x)\otimes 1)\Psi ({\bf 1})=Y(e^{x(1\otimes
\partial)}\Psi(v),x)({\bf 1}\otimes 1).$$
As $\Psi({\bf 1})=0$, we get
$$\Psi Y(v,x){\bf 1}=Y(e^{x(1\otimes \partial)}\Psi(v),x)({\bf 1}\otimes 1),$$
which implies
$$\Psi \D (v)=(\D\otimes 1)\Psi(v)+(1\otimes \partial)\Psi(v),$$
as we need.
\end{proof}

{}From the argument in the second part we also get:

\bl{ldproperty-end} Let $\Delta$ be an $R$-valued
pseudo-endomorphism of $V$. Then
\begin{eqnarray}
(\D\otimes 1)\Delta-\Delta \D=-(1\otimes
\partial)\Delta.
\end{eqnarray}
\el

As $R$ is a commutative associative algebra over $\C$ and $\partial$
is a derivation, by a result of Borcherds \cite{bor}, $R$ has a
vertex algebra structure (over $\C$) with $1$ as the vacuum vector
and with
\begin{eqnarray}
Y(a,x)b=(e^{x\partial}a)b=\sum_{n\ge
0}\frac{x^{n}}{n!}(\partial^{n}a)b
\end{eqnarray}
for $a,b\in R$. For convenience, we denote this vertex algebra by
$(R,\partial)$. Then for any vertex algebra $V$, we have a tensor
product vertex algebra $V\otimes (R,\partial)$, where the vertex
operator map, denoted by $Y_{ten}$, is given by
\begin{eqnarray}
Y_{ten}(u\otimes a,x)(v\otimes b)=Y(u,x)v\otimes
Y(a,x)b=Y(u,x)v\otimes (e^{x\partial}a)b
\end{eqnarray}
for $u,v\in V,\; a,b\in R$. Identify $V$ as a vertex subalgebra of
$V\otimes (R,\partial)$ in the canonical way, so that $V\otimes
(R,\partial)$ is a natural $V$-module. A linear map $\Phi:
V\rightarrow V\otimes R$ is called a {\em derivation} if
\begin{eqnarray}
\Phi Y(u,x)v=Y(u,x)\Psi (v)+Y_{ten}(\Phi (u),x)v
\end{eqnarray}
for $u,v\in V$.

Now, we interpret $R$-valued pseudo-derivations and -endomorphisms
as ordinary derivations and homomorphisms of vertex algebras.

\bl{lRpseudo} Let $V$ be a vertex algebra. An $R$-valued
pseudo-derivation (resp. pseudo-endomorphism) of $V$ exactly amounts
to a derivation (resp. an endomorphism) from $V$ to the tensor
product vertex algebra $V\otimes (R,\partial)$. \el

\begin{proof} Let $\Psi: V\rightarrow V\otimes R$ be a linear map. Notice that
$$Y(e^{x(1\otimes \partial)}(w\otimes a),x)=Y(w,x)\otimes Y(a,x)=Y_{ten}(w\otimes a,x)$$
for $w\in V,\; a\in R$. Then (\ref{e-pseudo-der}) exactly amounts to
$$ \Psi Y(v,x)u-(Y(v,x)\otimes 1)\Psi(u)=Y_{ten}(\Psi(v),x)(u\otimes 1).$$
Thus $\Psi$ is an $R$-valued pseudo-derivation if and only if $\Psi$
is a derivation of $V$-modules. The assertion for
pseudo-endomorphism is also clear.
\end{proof}

Using Lemma \ref{lRpseudo} we have (cf. \cite{li-pseudo}):

\bp{pordinary} Let $V$ be a vertex algebra. The restriction of every
derivation of the tensor product vertex algebra $V\otimes
(R,\partial)$ is an $R$-valued pseudo-derivation (resp.
pseudo-endomorphism) of $V$. Conversely, the unique $R$-linear
extension of any pseudo-derivation (resp. pseudo-endomorphism) of
$V$ is a derivation (resp. endomorphism) of $V\otimes (R,\partial)$.
\ep

\bd{dpsedo-auto}{\em A pseudo-endomorphism $\Delta$ of a vertex
algebra $V$ is called a {\em pseudo-automorphism} if the $R$-linear
extension of $\Delta$ is an automorphism of the tensor product
vertex algebra $V\otimes (R,\partial)$. } \ed

\br{rexample} {\em Consider the case $R=\C((z))$ with
$\partial=\frac{d}{dz}$. We denote a pseudo-derivation by $\Psi(z)$
to show its dependence of variable $z$. As
\begin{eqnarray*}
e^{x\frac{d}{d z}}\Psi(z)=\Psi(z+x),
\end{eqnarray*}
we have
\begin{eqnarray*}
\Psi(z) Y(v,x)-(Y(v,x)\otimes 1)\Psi(z)=Y(\Psi(z+x)v,x)
\end{eqnarray*}
for $v\in V$. With Lemma \ref{ldproperty-der}, we see that a
$\C((z))$-valued pseudo-derivation is exactly a pseudo-derivation in
the sense of \cite{ek}.  Similarly, a $\C((z))$-valued
pseudo-endomorphism is exactly a pseudo-endomorphism in the sense of
\cite{li-pseudo}.} \er

\br{rold} {\em Let $V$ be a vertex operator algebra and let $h\in
V_{(1)}$ be such that
$$L(n)h=\delta_{n,0}h,\ \ \ h_{n}h=\delta_{n,1}\gamma {\bf 1}$$
for $n\ge 0$, where $\gamma$ is a rational number. Assume that
$h_{0}$ acts on $V$ semisimply with eigenvalues in $\frac{1}{r}\Z$
for some positive integer $r$. Recall
$$\Delta(h,z)=z^{h(0)}\exp\left(\sum_{k=1}^{\infty}\frac{h(k)}{-k}
(-z)^{-k}\right)\in ({\rm End}V)[[z^{1/r},z^{-1/r}]].$$ The proof of
Proposition 5.4 of \cite{li-twisted} is mainly to prove that
$\Delta(h,z)$ is a $\C((z^{1/r}))$-valued pseudo-endomorphism.} \er

Just as with the classical case, the exponentials of
pseudo-derivations should give rise to pseudo-automorphisms.

\bl{lder-end}  Let $V$ be a vertex algebra and let $\Phi$ be an
$R$-valued pseudo-derivation of $V$, satisfying the condition that
for every $v\in V$, there exists a finite-dimensional subspace $U$
of $V$, containing $v$, such that
$$\Phi^{n+k}(U)\subset U\otimes I^{n}\ \ \ \mbox{ for }n\ge 0,$$
where $k$ is a fixed nonnegative integer and $I$ is a fixed ideal of
$R$. Suppose that $R$ is $I$-adically complete. Then
\begin{eqnarray}
\exp \Phi=\sum_{n\ge 0}\frac{1}{n!}\Phi^{n}
\end{eqnarray}
is an $R$-valued pseudo-automorphism of $V$. \el

\begin{proof} Notice that the assumptions imply that $\sum_{n\ge
0}\frac{1}{n!}\Phi^{n}(v)$ lies in $V\otimes R$ for $v\in V$, so
that $\exp \Phi$ is an $R$-linear automorphism of $V\otimes R$. By
Proposition \ref{pordinary}, $\Phi$, after $R$-linearly extended, is
a derivation of the vertex algebra $V\otimes (R,\partial)$. Since
$\Phi {\bf 1}=0$ (by Lemma \ref{ldproperty-der}), we have $(\exp
\Phi){\bf 1}={\bf 1}$. Let $m\in \Z$. Just as with any
non-associative algebra,  we have
$$(\exp \Phi) (a_{m}b)=(\exp \Phi)(a)_{m}(\exp \Phi)(b)
\ \ \ \mbox{ for }a,b\in V\otimes R.$$ This shows that $\exp \Phi$
is an automorphism of vertex algebra $V\otimes (R,\partial)$. By
Proposition \ref{pordinary}, $\exp \Phi$ is an $R$-valued
pseudo-automorphism of $V$.
\end{proof}

The following result, which slightly generalizes Proposition 1.9 of
\cite{ek} with a different proof, enables us to construct a family
of pseudo-derivations:

\bp{pconcrete} Let $V$ be a vertex algebra. For any $v\in V,\; f\in
R$, set
$$X_{v,f}=\sum_{n\ge 0}\frac{1}{n!}(\partial^{n}f)v_{n}.$$
Then $X_{v,f}$ is an $R$-valued pseudo-derivation of $V$. \ep

\begin{proof} Note that $(v\otimes f)_{0}$, which is the coefficient of $x^{-1}$
in the vertex operator $Y_{ten}(v\otimes f,x)$, is a derivation of
$V\otimes (R,\partial)$. In view of Proposition \ref{pordinary},
$(v\otimes f)_{0}$ is an $R$-valued pseudo-derivation of $V$. On the
other hand, we have
$$(v\otimes f)_{0}=\Res_{x}Y_{ten}(v\otimes
f,x)=\Res_{x}Y(v,x)\otimes e^{x\partial}f=X_{v,f}.$$ Now the
assertion follows.
\end{proof}

In view of Proposition \ref{pconcrete} and Lemma \ref{lder-end} we
immediately have:

\bc{cconstruction} Let $V$ be a vertex algebra and let $v\in V$.
Suppose that for any $w\in V$, there exists a positive integer $k$
such that
$$v_{n_{1}}v_{n_{2}}\cdots v_{n_{r}}w=0$$
for any nonnegative integers $n_{1},\dots,n_{r}$ with $n_{1}+\cdots
+n_{r}\ge k$. Then for any $f\in R$,
$$\exp \left(\sum_{n\ge 0}\frac{1}{n!}(\partial^{n}f)v_{n}\right)$$
is an $R$-valued pseudo-automorphism of $V$. In particular, this is
true if $V$ is a vertex operator algebra and if $v$ is homogeneous
with a non-positive conformal weight. \ec

\begin{proof} The first assertion is immediate. For the second assertion,
as $\wt\; v\le 0$, we have $\wt\; v_{n}=\wt\; v-n-1<0$ for $n\ge 0$.
Since $V$ is truncated from below by definition, the assumption in
the first part is satisfied. Then it follows.
\end{proof}

\bp{pexponential} Let $V$ be a vertex operator algebra and let $v\in
V$ be homogeneous of weight $k$ positive, such that
$[v_{m},v_{n}]=0$ for $m,n\ge 0$. Suppose that
$v_{0},v_{1},\dots,v_{k-1}$ are locally finite on $V$ and suppose
that $f(z)\in z^{k}\C[[z]]$. Then $\exp X_{v,f}$ exists in $\Hom
(V,V\otimes \C((z)))$ and it is a $\C((z))$-valued
pseudo-automorphism of $V$. \ep

\begin{proof} As $f^{(n)}(z)\in z\C[[z]]$ for $0\le n\le k-1$ and as
$v_{0},v_{1},\dots,v_{k-1}$ are locally finite, we see that
$$\exp \left(\sum_{n=0}^{k-1}\frac{1}{n!}f^{(n)}(z)v_{n}\right)$$
exists in $\Hom (V,V\otimes \C[[z]])$ (not just in $(\End V)[[z]]$).
On the other hand,
 $$\exp \left(\sum_{n\ge k}\frac{1}{n!}f^{(n)}(z)v_{n}\right)$$
also exists in $\Hom (V,V\otimes \C[[z]])$ as $\wt\;v_{n}=k-n-1<0$
for $n\ge k$. Consequently, $\exp X_{v,f}$ exists in $\Hom
(V,V\otimes \C((z)))$ and it is a $\C((z))$-valued
pseudo-automorphism.
\end{proof}

\bex{example-fraction} {\em Let $V$ be a vertex operator algebra and
let $v\in V_{(1)}$ be such that $[v_{m},v_{n}]=0$ for $m,n\ge 0$ and
such that $v_{0}$ is locally finite on $V$. Let $r$ be a positive
integer and let $f(z)=z^{1/r}$. We have
$$X_{v,f}=\sum_{n\ge 0}\binom{\frac{1}{r}}{n}v_{n}z^{-n+\frac{1}{r}}.$$
Note that $\exp (v_{0}z^{1/r})$ exists in $\Hom (V,V\otimes
\C[[z^{1/r}]])$ and $\exp\left(\sum_{n\ge
1}\binom{\frac{1}{r}}{n}v_{n}z^{-n+\frac{1}{r}}\right)$ exists in
$\Hom (V,V\otimes \C[z^{-1/r}])$. Then $\exp (X_{v,f})$ exists in
$\Hom (V,V\otimes \C((z^{1/r})))$ and it is a $\C((z^{1/r}))$-valued
pseudo-automorphism of $V$.} \eex

The importance of pseudo-endomorphisms can be seen from the
following results due to \cite{li-selection} and \cite{dlm-simple}
(cf. \cite{li-pseudo}):

\bp{pmodule-sele} Let $V$ be a vertex algebra and let $\Delta(z)$ be
a ($\C((z))$-valued) pseudo-endomorphism of $V$. Then for any
$V$-module $(W,Y_{W})$, $(W,Y_{W}^{\Delta})$ carries the structure
of a $V$-module where
$$Y_{W}^{\Delta}(v,x)=Y_{W}(\Delta(x)v,x)\ \ \ \mbox{ for }v\in V.$$
 \ep

\bp{ptwistedmodule-sele} Let $V$ be a vertex algebra and let
$\sigma$ be an automorphism of $V$ of order $r$. Suppose that
$\Delta(z)$ is a $\C((z^{1/r}))$-valued pseudo-endomorphism of $V$,
satisfying
\begin{eqnarray}\label{esigma-tau}
(\sigma\times \tau)\Delta(z)=\Delta(z),
\end{eqnarray}
where $\tau f(z^{1/r})=f(\omega_{r}z^{1/r})$ for $f(z)\in \C((z))$.
Then for any $V$-module $(W,Y_{W})$, $(W,Y_{W}^{\Delta})$ carries
the structure of a $\sigma$-twisted $V$-module where
$$Y_{W}^{\Delta}(v,x)=Y_{W}(\Delta(x)v,x)\ \ \ \mbox{ for }v\in V.$$
 \ep

\bex{M(1)} {\em  Let $V=M(1)$ be the (free field) Heisenberg vertex
operator algebra with one generator $h$ satisfying
$$[h_{m},h_{n}]=m\delta_{m+n,0}$$
for $m,n\in \Z$. Note that $h_{0}=0$ on $M(1)$. We see that the
assumption in Corollary \ref{cconstruction} with $v=h$ is satisfied.
Thus for any $f\in \C((z))$, $\exp X_{h,f}$, lying in $\Hom
(M(1),M(1)\otimes \C((z)))$, is a pseudo-automorphism of $M(1)$. Let
$k$ be a positive integer and set
$$\Phi_{k}(z)=\exp
\left(X_{h,z^{-k}}\right)=\exp \left(\sum_{n\ge
1}\binom{-k}{n}h(n)z^{-n-k}\right).$$  Noticing that
$X_{h,z^{-k}}(h)=-kz^{-k-1}{\bf 1}$, we have
$\Phi_{k}(z)(h)=h-kz^{-k-1}{\bf 1}$. For any $M(1)$-module
$(W,Y_{W})$, we have a module $(W,Y_{W}^{\Phi_{k}})$ for $M(1)$
viewed as a vertex algebra, where
$$Y_{W}^{\Phi_{k}}(h,z)=Y_{W}(\Phi_{k}(z)h,z)=Y_{W}(h,z)-kz^{-k-1}.$$
Set $Y_{W}^{\Phi_{k}}(h,z)=\sum_{n\in \Z}h^{\Phi_{k}}_{n}z^{-n-1}$.
We have
$$h^{\Phi_{k}}_{n}=h_{n}-k\delta_{n,k}\ \ \ \mbox{ for }n\in \Z.$$
Take $W=M(1)$, the adjoint module, on which $h_{k}$ is locally
nilpotent. We see that $h^{\Phi_{k}}_{k}=h_{k}-k$ is locally finite
but not nilpotent. Thus the deformed $M(1)$-module
$(M(1),Y^{\Phi_{k}})$ is not a highest weight module. Such
$M(1)$-modules viewed as modules for the Heisenberg algebra were
studied in \cite{lw}.} \eex

Next, we interpret $\Delta(h,z)$ (recall Remark \ref{rold}) as the
exponential of a pseudo-derivation. First we define an algebra $R$
involving $\log z$. Let $r$ be a positive integer. Consider algebra
$\C((z^{1/r}))[[y]]$ with $y$ an independent variable. There exists
a $\C$-linear derivation $\partial$ such that
\begin{eqnarray}
\partial\left(\sum_{n\ge 0}g_{n}y^{n}\right)=\sum_{n\ge
0}\left(\frac{d g_{n}}{dz}y^{n}+nz^{-1}g_{n}y^{n-1}\right)
\end{eqnarray}
for $\sum_{n\ge 0}g_{n}y^{n}\in \C((z^{1/r}))[[y]]$. In particular,
we have $\partial y=z^{-1}$. The algebra $\C((z^{1/r}))[[y]]$
contains $\C((z))\otimes \C[[y]]$, $\C[z,z^{-1}]\otimes \C[[y]]$,
and $\C[y, z^{-1}]$ as subalgebras which are all stable under the
derivation operator $\partial$.

Set
\begin{eqnarray}
K_{r}[y]=\left\<y,\ e^{\alpha y}\;|\; \alpha \in
\frac{1}{r}\Z\right>,
\end{eqnarray}
a subalgebra of $\C[[y]]$, where
$$e^{\alpha y}=\sum_{n\ge 0}\frac{1}{n!}\alpha^{n}z^{n}\in
\C[[z]].$$ One can show that $K_{r}[y]$ is isomorphic to
$\C[\frac{1}{r}\Z]\otimes \C[y]$, where $\C[\frac{1}{r}\Z]$ stands
for the group algebra of $\frac{1}{r}\Z$.

Note that both $\C((z))\otimes K_{r}[y]$ and $\C((z^{1/r}))[y]$ are
subalgebras of $\C((z^{1/r}))[[y]]$, which are stable under
$\partial$. The following is straightforward:

\bl{lalgebra-1} There exists an algebra homomorphism
$$\theta: \C((z))\otimes K_{r}[y]\rightarrow \C((z^{1/r}))[y],$$
which is uniquely determined by
$$\theta(f(z))=f(z),\ \ \ \theta(y)=y,\ \ \ \theta( e^{\pm y/r})=z^{\pm 1/r}$$
for $f(z)\in \C((z))$. Furthermore, we have  $\theta
\partial=\partial \theta$.
 \el

For any vertex algebra $V$ and for any $v\in V$, by Proposition
\ref{pconcrete},
$$X_{v,y}=\sum_{n\ge 0}\frac{1}{n!}(\partial^{n}y)
v_{n}=v_{0}y+\sum_{n\ge 1}(-1)^{n-1}\frac{1}{n}v_{n}z^{-n}$$ is a
$\C[y,z^{-1}]$-valued pseudo-derivation of $V$.

\bp{palgebra-0} Let $V$ be a vertex operator algebra and let $h\in
V_{(1)}$, satisfying that $[h_{m},h_{n}]=0$ for $m,n\ge 0$ and that
$h_{0}$ is locally finite on $V$. Then $\exp (X_{h,y})$ is a
$\C[z^{-1}]\otimes \C[[y]]$-valued pseudo-automorphism of $V$.
Furthermore, if the eigenvalues of $h_{0}$ lie in $\frac{1}{r}\Z$,
$\exp (X_{h,z})$ is a $\C[z^{-1}]\otimes K_{r}[y]$-valued
pseudo-endomorphism of $V$. \ep

\begin{proof} Since $\wt\; h_{n}=-n\le 1$ for $n\ge 1$, we see that $\sum_{n\ge
1}(-1)^{n-1}\frac{1}{n}h_{n}z^{-n}$ is locally nilpotent on $V$. As
$h_{0}$ is assumed to be locally finite on $V$, we have
$$e^{yh_{0}}\in \Hom (V,V\otimes \C[[y]]).$$
Thus
$$\exp (X_{h,y})\in \Hom (V,V\otimes \C[z^{-1}]\otimes \C[[y]]).$$
 We see that $\exp (X_{h,y})$ is a $\C[z^{-1}]\otimes \C[[y]]$-valued
pseudo-endomorphism of $V$.

Now, assume that the eigenvalues of $h_{0}$ lie in $\frac{1}{r}\Z$.
By using the semisimple-nilpotent decomposition of $h_{0}$, we have
$$e^{yh_{0}}\in \Hom (V,V\otimes K_{r}[y]).$$
Consequently, we obtain
$$\exp (X_{h,y})\in \Hom (V,V\otimes \C[z^{-1}]\otimes K_{r}[y]),$$
proving that $\exp (X_{h,z})$ is a $\C[z^{-1}]\otimes
K_{r}[y]$-valued pseudo-endomorphism of $V$.
\end{proof}

Identify $\C((z^{1/r}))[\log z]$ with the subalgebra
$\C((z^{1/r}))[y]$ of $\C((z^{1/r}))[[y]]$ where $y=\log z$. That
is, $\log z$ is considered as a formal variable. We have (cf.
\cite{am}, \cite{huang}):

\bc{chuang} Let $V$ be a vertex operator algebra and let $h\in
V_{(1)}$, satisfying that $[h_{m},h_{n}]=0$ for $m,n\ge 0$ and that
$h_{0}$ is locally finite on $V$ with eigenvalues lying in
$\frac{1}{r}\Z$. Set
$$\Delta(h,z)=\exp(h_{0}\log z)\exp \left(\sum_{n\ge
1}(-1)^{n-1}\frac{1}{n}h_{n}z^{-n} \right),$$ where for $v\in V$,
$e^{h_{0}\log z}v$ is defined by
$$e^{h_{0}\log z}v=z^{h_{0}^{ss}}e^{h_{0}^{n}}v,  $$
where  $h_{0}=h_{0}^{ss}+h_{0}^{n}$ is the semisimple-nilpotent
decomposition on any finite-dimensional $h_{0}$-stable subspace $U$
containing $v$. Then  $\Delta(h,z)$ is a $\C((z^{1/r}))[\log
z]$-valued pseudo-endomorphism of $V$. \ec

\begin{proof} It is clear that the definition does not depend on the choice of $U$.
By Proposition \ref{palgebra-0}, $\exp (X_{h,y})$ is a
$\C[z^{-1}]\otimes K_{r}[y]$-valued pseudo-endomorphism of $V$. In
view of Lemma \ref{lRpseudo}, $\exp (X_{h,y})$ is a vertex algebra
homomorphism from $V$ to $V\otimes (\C[z^{-1}]\otimes
K_{r}[y],\partial)$. On the other hand, {}from Lemma
\ref{lalgebra-1}, we see that the algebra homomorphism $\theta$ is
actually a vertex algebra homomorphism {}from $(\C((z))\otimes
K_{r}[y],\partial)$ to $(\C((z^{1/r}))[y],\partial)$. Then $\theta
\circ \exp (X_{h,y})$ is a vertex algebra homomorphism from $V$ to
$V\otimes (\C((z^{1/r}))[y],\partial)$ and by Lemma \ref{lRpseudo}
it is a $\C((z^{1/r}))[y]$-valued pseudo-endomorphism of $V$.
\end{proof}

Now we consider pseudo-endomorphism twists of representations for a
vertex algebra. Recall from \cite{ltw} (Remark 2.16) the following
result on the definition of a twisted module (cf.
\cite{li-twisted}):

\bp{pltw} Let $V$ be a vertex algebra and let $\sigma$ be an
automorphism of $V$ of order $r$. Set
$$\omega_{r}=\exp (2\pi \sqrt{-1}/r).$$
Let $W$ be a vector space and let
$$Y_{W}: V\rightarrow \Hom (W,W((x^{\frac{1}{r}})))
\subset (\End W)[[x^{\frac{1}{r}},x^{-\frac{1}{r}}]]$$ be a linear
map. Then $(W,Y_{W})$ carries the structure of a  $\sigma$-twisted
$V$-module if and only if the following conditions hold: $Y_{W}({\bf
1},x)=1_{W}$,
\begin{eqnarray}\label{esigma-invariance}
Y_{W}(\sigma v,x)=\lim_{x^{\frac{1}{r}}\rightarrow
\omega_{r}^{-1}x^{\frac{1}{r}}}Y_{W}(v,x)
\end{eqnarray}
 for $v\in V$, and for $u,v\in V$, there exists $k\in \N$ such that
 $$(x_{1}-x_{2})^{k}Y_{W}(u,x_{1})Y_{W}(v,x_{2})\in
\Hom (W,W((x_{1}^{\frac{1}{r}},x_{2}^{\frac{1}{r}})))$$ and
\begin{eqnarray*}
x_{0}^{k}Y_{W}(Y(u,x_{0})v,x_{2})
=\left((x_{1}-x_{2})^{k}Y_{W}(u,x_{1})Y_{W}(v,x_{2})\right)|_{x_{1}^{1/r}=(x_{2}+x_{0})^{1/r}}.
\end{eqnarray*}
\ep

\bd{dr-twisted} {\em Let $V$ be a vertex algebra and let $r$ be a
positive integer. A {\em $\C((z^{1/r}))$-valued $V$-module} is a
vector space $W$ equipped with a linear map
 $$Y_{W}: V\rightarrow \Hom (W,W((x^{1/r})))
 \subset (\End W)[[x^{\frac{1}{r}},x^{-\frac{1}{r}}]],$$
satisfying all the conditions listed in Proposition \ref{pltw}
except (\ref{esigma-invariance}). } \ed

The following generalizes Proposition \ref{ptwistedmodule-sele} in a
certain way (cf. \cite{huang}):

\bp{ptwistedmodule-deformation} Let $V$ be a vertex algebra and let
$\Delta(z)$ be a $\C((z^{1/r}))$-valued pseudo-endomorphism of $V$
with $r$ a positive integer. Let $(W,Y_{W})$ be a $V$-module. For
$v\in V$, define
\begin{eqnarray}
Y_{W}^{\Delta}(v,x)=Y_{W}(\Delta(x)v,x).
\end{eqnarray}
Then $(W,Y_{W}^{\Delta})$ carries the structure of a
$\C((z^{1/r}))$-valued $V$-module.\ep

\begin{proof} For $v\in V,\; w\in W$, as $\Delta(x)v\in V\otimes \C((x^{1/r}))$, we
have $Y_{W}^{\Delta}(v,x)w\in W((x^{1/r}))$.
 We have
$$Y_{W}^{\Delta}({\bf 1},x)=Y_{W}(\Delta(x)({\bf 1}),x)=Y_{W}({\bf
1},x)=1_{W}.$$ Let $u,v\in V$. As $\Delta(z)(u),\Delta(z)(v)\in
V\otimes \C((z^{1/r}))$, there exists a nonnegative integer $k$ such
that
$$(x_{1}-x_{2})^{k}[Y_{W}(\Delta(z_{1})u,x_{1}),Y_{W}(\Delta(z_{2})v,x_{2})]=0.$$
Substituting $z_{1}^{1/r}=x_{1}^{1/r},\; z_{2}^{1/r}=x_{2}^{1/r}$,
we get
$$(x_{1}-x_{2})^{k}[Y_{W}^{\Delta}(u,x_{1}),Y_{W}^{\Delta}(v,x_{2})]=0,$$
which implies
$$(x_{1}-x_{2})^{k}Y_{W}^{\Delta}(u,x_{1})Y_{W}^{\Delta}(v,x_{2})\in
\Hom (W,W((x_{1}^{1/r},x_{2}^{1/r}))).$$ We also have
\begin{eqnarray*}
&&x_{0}^{k}Y_{W}(Y(\Delta(z_{1})u,x_{0})\Delta(z_{2})v,x_{2}) \\
&=&\left((x_{1}-x_{2})^{k}Y_{W}(\Delta(z_{1})u,x_{1})
Y_{W}(\Delta(z_{2})v,x_{2})\right)|_{x_{1}^{1/r}=(x_{2}+x_{0})^{1/r}}.
\end{eqnarray*}
Substituting $z_{1}^{1/r}=(x_{2}+x_{0})^{1/r},\;
z_{2}^{1/r}=x_{2}^{1/r}$, we get
\begin{eqnarray*}
x_{0}^{k}Y_{W}^{\Delta}(Y(u,x_{0})v,x_{2})
=\left((x_{1}-x_{2})^{k}Y_{W}^{\Delta}(u,x_{1})
Y_{W}^{\Delta}(v,x_{2})\right)|_{x_{1}^{1/r}=(x_{2}+x_{0})^{1/r}},
\end{eqnarray*}
noticing that
\begin{eqnarray*}
&&x_{0}^{k}Y_{W}(Y(\Delta(x_{2}+x_{0})u,x_{0})\Delta(x_{2})v,x_{2})\\
&=&x_{0}^{k}Y_{W}(\Delta(x_{2})Y(u,x_{0})v,x_{2})\\
&=&x_{0}^{k}Y_{W}^{\Delta}(Y(u,x_{0})v,x_{2}).
\end{eqnarray*}
This completes the proof.
\end{proof}

\br{rR-valuedmodule} {\em Let $R$ be a commutative algebra over
$\C((z))$ as before. For a vertex algebra $V$, one can define a
notion of $R$-valued $V$-module. } \er

\end{document}